\documentclass[a4paper,11pt]{article}
\usepackage[utf8]{inputenc}
\usepackage{color}

\usepackage[all]{xy}
\usepackage{amsmath,amssymb,amsfonts,amsthm,graphicx}
\usepackage[margin=1.15in]{geometry}
\usepackage{multirow}
\theoremstyle{plain}

\newtheorem*{theorem*}{Theorem}
\newtheorem*{aim*}{Aim}
\newtheorem*{initialaim*}{Initial Aim}
\newtheorem*{conj*}{Conjecture}
\newtheorem*{cor*}{Corollary}
\newtheorem*{prop*}{Proposition}
\newtheorem*{df*}{Definition}
\newtheorem*{lm*}{Lemma}
\newtheorem*{example*}{Example}
\newtheorem*{notation*}{Notation}
\newtheorem*{prob*}{Problem}

\newcommand{\Z}{\mathbb{Z}}

\numberwithin{thm}{subsection}

\title{Non-vanishing elements in finite groups}
\author{Julian Brough}

\begin{document}
\date{}
\maketitle

\begin{center}
\small
\textit{FB Mathematik, TU Kaiserslautern, Postfach 3049}

\textit{67653 Kaiserslautern, Germany}
\end{center}

\paragraph{}

\normalsize
\begin{abstract}
Many results have been established about determining whether or not an element evaluates to zero on an irreducible character of a group.
In this note it is shown that if a group $G$ has a normal nilpotent subgroup $N$, and $P$ is a Sylow $p$-subgroup of $G$, then no irreducible character of $G$ vanishes on $N\cap Z(P)$.
\end{abstract}
Let $G$ be a finite group and $\chi\in$ Irr$(G)$, an irreducible character of $G$.
A classical result of Burnside says if $\chi$ is non-linear, that is $\chi(1)\not=1$, then there is at least one element $g$ in $G$ such that $\chi(g)=0$.
If one considers conjugacy classes, a natural dual to irreducible characters, then $g$ being a central element in $G$ implies that $|\chi(g)|=\chi(1)$ and thus $g$ does not evaluate to zero on any irreducible character.
However, a non-central element $g$ may also not evaluate to zero on any irreducible character, for example the $3$-cycles in ${\rm Sym}(3)$.
Elements which do not evaluate to zero on any irreducible character of a group are called non-vanishing.
The study of non-vanishing elements was first introduced in \cite{INW}, where the authors showed for soluble groups any non-vanishing element $g$ in a group $G$ must reduce to a $2$-element in $G/F(G)$. 
In \cite{DNPST} this result was generalised to any group, in particular it was shown that if an element $g$ is non-vanishing in $G$ and the order of $g$ is coprime to $6$, then $g$ lies in $F(G)$.

We note that there has been a recent interest in the literature asking about how much group structure is determined by the vanishing conjugacy class sizes.
In particular, in \cite{DPSVan} and \cite{Brough3}, the authors have generalised arithmetical results upon conjugacy classes to vanishing conjugacy classes.
Thus the determination of non-vanishing elements would provide further machinery for this recent topic of research.
 
The aim of this note is to generalise one of the key results in \cite{INW}, that is \cite[Theorem A]{INW}, which says if a group has a normal Sylow $p$-subgroup $P$, then all the elements in $Z(P)$ are non-vanishing.
A variant of this result was considered in \cite{NonVanMas}, where the author showed that if a group has a normal elementary abelian $p$-subgroup $A$ and $P$ is a Sylow $p$-subgroup, then the elements in $Z(P)\cap A$ are non-vanishing.
In particular, we first show that the result of \cite{NonVanMas} holds if $A$ is a normal abelian subgroup.
From this we deduce the result holds if $A$ is a normal nilpotent subgroup.
Note that from this new result, \cite[Theorem A]{INW} follows by setting $A=P$.

\begin{theorem*}
Let $G$ be a finite group, which contains a non-trivial normal nilpotent subgroup $N$ and $p$ a prime.
Then for $P\in Syl_p(G)$, the elements in $N\cap Z(P)$ are non-vanishing in $G$.
\end{theorem*}

First we give the following preliminary result which considers when a sum of roots of unity is equal to zero.
\begin{lm*}
Let $\Xi:=\{\xi_1,\dots, \xi_t\}$ be a set of $p^n$-th roots of unity, for some number $n\geq 1$, such that $\xi_1+\dots+\xi_t=0$.
Then the sum can be split into sums of the form $\xi+\xi^{p^{a-1}+1}+\dots +\xi^{(p-1)p^{a-1}+1}$, for possibly various numbers $1\leq a\leq n$, where $\xi^{kp^{a-1}+1}\in\Xi$ for $0\leq k\leq p-1$ and each such subsum equals zero.
\begin{proof}
Let $\xi$ be an element in $\{\xi_i \mid 1\leq i\leq t\}$ of maximal order, i.e. $\xi$ is a primitive $p^a$-th root of unity and $\xi_i^{p^a}=1$ for all $i$.
It is enough to prove that $\xi,\xi^{p^{a-1}+1},\dots, \xi^{(p-1)p^{a-1}+1}\in \Xi$, as then
\[
\xi+\xi^{p^{a-1}+1}+\dots +\xi^{(p-1)p^{a-1}+1}=\xi(1+\xi^{p^{a-1}}+\dots +\xi^{(p-1)p^{a-1}})=0
\]
and inductively from $\Xi\backslash \{\xi,\xi^{p^{a-1}+1},\dots \xi^{(p-1)p^{a-1}+1}\}$ repeat the argument.

Assume $\xi_1=\xi$, which is a primitive $p^a$-th root of unity, so that each $\xi_i$ is a power of $\xi$.
Pick $r$ minimal such that $\Sigma_{j=1}^r\xi^{b_j}=0$ with $\xi^{b_j}\in \Xi$, where $b_{i}\leq b_{i+1}\leq p^a$ and $b_1=1$.
Then it follows that $\xi$ is a root to the polynomial $\Sigma_{j=1}^rX^{b_j}$. 
As $\Phi_{p^a}(X)=1+X^{p^{a-1}}+\dots +X^{(p-1)p^{a-1}}$ is the minimal polynomial for $\xi$ it follows that 
\[
\Sigma_{j=1}^rX^{b_j}=\Phi_{p^a}(X) g(X),
\]
for some polynomial $g\in \Z[X]$.

For each $j$, $b_j>0$, therefore $g(X)$ cannot have a constant term, i.e. $g(X)=Xf(X)$ for some $f\in \Z[X]$.
The polynomial $f(X)$ must have a constant term $c\ne 0$, because $b_1=1$.
Thus
\[
\Sigma_{j=1}^rX^{b_j}=cX\Phi_{p^a}(X)+ X^2\Phi_{p^a}(X)h(X),
\]
for some $h\in\Z[X]$.
As $b_j\leq p^a$, it follows that $2+(p-1)p^{a-1}+{\rm deg}(h(X))\leq p^a$.
Moreover, $2\leq {\rm deg}(X^2h(X))\leq p^{a-1}$.
If $X^2\Phi_{p^a}(X)h(X)$ has a term of the form $X^{kp^{a-1}+1}$ then $X^2h(X)$ must have the term $X$ or $X^{mp^{a-1}+1}$ for some positive integer $m$, which is a contradiction.
In particular, $cX\Phi_{p^a}(X)$ has no terms in common with $X^2\Phi_{p^a}(X)h(X)$.
Hence $X\Phi_{p^a}(X)$ occurs as a subsum of $\Sigma_{j=1}^rX^{b_j}$.
However as $\xi\Phi_{p^a}(\xi)=0$ it follows that
\[
\Sigma_{j=1}^rX^{b_j}=X\Phi_{p^a}(X).
\]
\end{proof}
\end{lm*}

We can now establish the main result in the case that a group has a normal abelian subgroup. 
Note that the proof makes use of the method in \cite{NonVanMas} with the additional information about roots of unity in the previous lemma.

\begin{prop*}
Let $G$ be a finite group, which contains a non-trivial normal abelian subgroup $A$ and $p$ a prime.
Then for $P\in Syl_p(G)$, the elements in $A\cap Z(P)$ are non-vanishing in $G$.
\begin{proof}
Let $x\in A\cap Z(P)$ such that there exists some $\chi\in {\rm Irr}(G)$ for which $\chi(x)=0$.
By Clifford's theorem $\chi\downarrow_A=e\Sigma_{i=1}^t\zeta^{g_i}$ such that $\zeta\in{\rm Irr}(A)$ and the set $\{g_i\}$ forms a transversal of $I:=I_G(\zeta)$ in $G$, for $I_G(\zeta)$ the inertial subgroup in $G$ of $\zeta$.
If $\chi\downarrow_A(x)=0$, then $\Sigma_{i=1}^t\zeta^{g_i}(x)=0$.
Thus by the lemma we can split this sum into smaller subsums of $p$ elements (which also equal zero).
Let $\{\xi_j\}$ denote a set of representatives for the distinct subsums of $p$ elements as in the above lemma.
If $k_j$ denotes the multiplicity of $\xi_j$, then $\Sigma_jk_jp=t$.
Hence it is enough to show that the $p$-part of $t$, denoted $t_p$, divides $k_j$ as then $t_pp$ divides $t$ which is a contradiction.

The fact that the multiplicity $k_j$ is divisible by $t_p$ is in the proof of \cite[Theorem]{NonVanMas}, however we shall include details for completeness.

The subgroup $P$ acts on the set of $G$-conjugates of $\zeta$ with the orbit size of $\zeta^g$ given by
\[
|P:P\cap I^g|=|G:P\cap I^g|_p=|G:I^g|_p|I^g:P\cap I^g|_p.
\]
Therefore $|G:I|_p$ divides the orbit size.
As $x\in Z(P)$ it follows that $\zeta^{gy}(x)=\zeta^g(x)$ for all $y\in P$.
Thus the value $x$ evaluates to is constant on each $P$-orbit.
Hence the multiplicity of $\xi$ in $\{\zeta^{g_i}(x)\}$ must be divisible by $|G:I|_p=t_p$.

This completes the proof.
\end{proof}
\end{prop*}

From the proposition the main theorem now follows.

\begin{proof}[{\bf Proof of theorem}]
Let $N$ be a non-trivial normal nilpotent subgroup of a finite group $G$ and $P$ a Sylow $p$-subgroup of $G$ for a prime $p$.
Assume $x\in N\cap Z(P)$.
As $N$ is nilpotent, $O_p(N)$ is the unique Sylow $p$-subgroup of $N$.
Thus $x\in O_p(N)\cap Z(P)$.
As $O_p(N)$ is a normal $p$-subgroup of $G$, the group $O_p(N)$ is a subgroup of $P$.
In particular, it follows that $O_p(N)\cap Z(P)\leq Z(O_p(N))$.
Hence $x\in Z(O_p(N))\cap Z(P)$ and the result follows from the proposition, as $Z(O_p(N))$ is a normal abelian subgroup of $G$.
\end{proof}

\begin{cor*}[{\bf 1}]
Let $G$ be a finite group and $F(G)$ the fitting subgroup of $G$.
Then for $p$ a prime and $P\in Syl_p(G)$, the elements in $F(G)\cap Z(P)$ are non-vanishing in $G$.
\end{cor*}

In \cite{INW} it is conjectured that for soluble groups, all non-vanishing elements lie in the Fitting group.
This result provides a partical insight into which elements of the Fitting subgroup would in fact be non-vanishing in a group. 

In the proof of the proposition, the importance of $A$ being abelian is that the restricted character must be a sum of linear characters. 
It is therefore natural to assume that if such a restriction was made to characters with degree not divisible by $p$, a similar argument should work.
In fact this can be bypassed in a more general setting by using the Ito-Michler Theorem, which was proven for soluble groups by N. Ito \cite{ItoGpChar} and then for any group by G. Michler \cite{MichBraCFSG} using the classification of finite simple groups.

\begin{cor*}[{\bf 2}]
Let $G$ be a finite group and $N$ normal in $G$ such that no irreducible character of $N$ has degree divisible by a prime $p$.
Furthermore, let $P$ be a Sylow $p$-subgroup of $G$.
Then the elements in $Z(P)\cap N$ are non-vanishing in $G$.
\begin{proof}
As $N$ has no irreducible characters of degree divisible by $p$, the Ito-Michler theorem implies that $N$ has a normal Sylow $p$-subgroup $Q$ which is abelian.
Therefore $N\cap Z(P)=Q\cap Z(P)$.
Moreover, as $Q$ is a normal Sylow $p$-subgroup of $N$, $Q$ is normal in $G$.
Hence $Q$ is a normal abelian subgroup of $G$.
Thus by the theorem, the elements of $Q\cap Z(P)$ are non-vanishing in $G$. 
\end{proof}
\end{cor*}

Finally note that in the above corollary the condition that no irreducible character has degree divisible by $p$ cannot be removed.
In particular, let $G$ be a non-abelian simple group of Lie type with order divisible by $p$.
Then $G$ has an irreducible character of $r$-defect zero for any prime $r$ and hence any non-trivial element in $G$ is vanishing \cite[Theorem 8.17]{IsaChTh}.
Thus for $G=N$ the conclusion of the corollary cannot hold.

\vspace{4mm}

\noindent{\bf Acknowledgements}

\vspace{3mm}

The author wishes to thank Professor Malle for a helpful discussion during the preparation of this paper.  

\bibliographystyle{alpha}
\bibliography{bibfile1}

\end{document}